\documentclass{amsart}

\newtheorem{theorem}{Theorem}
\newcommand{\bt}{\begin{theorem}}
\newcommand{\et}{\end{theorem}}

\newtheorem{problem}{Problem}
\newcommand{\bprob}{\begin{problem}}
\newcommand{\eprob}{\end{problem}}

\newcommand{\mcr}{\ensuremath{ \mathcal R}} 

\newcommand{\N}{\ensuremath{ \mathbf N }}
\newcommand{\Z}{\ensuremath{\mathbf Z}}

\newcommand{\R}{\ensuremath{\mathbf R}}

\newcommand{\mbx}{\ensuremath{\mathbf x}}
\newcommand{\mby}{\ensuremath{\mathbf y}}
\newcommand{\mba}{\ensuremath{\mathbf a}}
\newcommand{\mbb}{\ensuremath{\mathbf b}}

\newcommand{\beq}{\begin{equation}}
\newcommand{\eeq}{\end{equation}}

\DeclareMathOperator{\qqand}{\qquad\text{and}\qquad}

\DeclareMathOperator{\diam}{\text{diam}}

\title{A problem on sumset sizes of sets of lattice points}

\author{Melvyn B. Nathanson}

\address{Department of Mathematics, Lehman College (CUNY), Bronx, NY, USA}  
\email{melvyn.nathanson@lehman.cuny.edu}

\date{\today}

\begin{document}

\maketitle

\begin{abstract}
A central problem in additive number theory is to understand the set of sizes of $h$-fold sums 
of finite subsets of an additive abelian semigroup.  It is proved that these ``range of sumset sizes'' 
is the same for finite sets of integers and for finite sets of $n$-dimensional lattice points.  
The associated geometrical problem is to determine if these sumset size sets 
can be computed more efficiently 
with lattice points than with integers. 
\end{abstract}

Let $A$ be a nonempty subset of an additive abelian group or semigroup $X$ 
and let $h$ and $k$ be positive integers. 
The \emph{$h$-fold sumset of $A$} is the set $hA$ of all sums of $h$ not necessarily 
distinct elements of $A$. 
The \emph{restricted $h$-fold sumset of $A$} is the set $hA$ of all sums of $h$ 
pairwise distinct elements of $A$. 

Nathanson~\cite{nath25a} introduced the study of sumset sizes of finite sets.  
For all positive integers $h$ and $k$, we define the 
 the \emph{range of sumset sizes} of finite subsets 
of  the semigroup $X$: 
\[
\mcr_{X} (h,k) = \left\{ \left| hA \right|: A \subseteq X \text{ and }  |A| = k \right\} 
 \]
and  the \emph{range of restricted sumset sizes} of finite subsets of $X$:  
\[
\widehat{\mcr}_X (h,k) = \left\{ \left|\widehat{hA} \right|: A \subseteq X \text{ and }  |A| = k \right\}.
 \] 
If $X$ is a isomorphic to a subsemigroup of the semigroup $Y$, then 
\[
\mcr_{X} (h,k) \subseteq \mcr_{Y} (h,k) 
 \]
and 
\[
\widehat{\mcr}_X (h,k) \subseteq \widehat{\mcr}_Y (h,k).
 \] 
It is of interest to consider  pairs of semigroups $X$ and $Y$ such that 
\[
\mcr_{X} (h,k) = \mcr_{Y} (h,k) 
 \]
and 
\[
\widehat{\mcr}_X (h,k) = \widehat{\mcr}_Y (h,k).
 \] 
for some or all pairs $(h,k)$. 
For example,  the additive group $\Z$ of integers is isomorphic 
to a subgroup of the the group $\Z^n$ of $n$-dimensional lattice points and so 
\[
\mcr_{\Z} (h,k) \subseteq \mcr_{\Z^n} (h,k).
 \]

\bt           \label{lattice:theorem:Z}
For all positive integers $h$ and $k$, 
\[
\mcr_{\Z} (h,k) = \mcr_{\Z^n} (h,k) 
 \] 
 and 
 \[
\widehat{\mcr}_{\Z} (h,k) = \widehat{\mcr}_{\Z^n} (h,k).
 \]
\et 

\begin{proof}
Let $g$, $h$, and $\ell$ be positive integers with 
\[
g > h\ell.  
\]
We define the linear function   $\varphi_g:\R^n \rightarrow \R$  by 
\[
\varphi_g(x_1,\ldots, x_n) = \sum_{i=1}^n x_i g^{i-1}.
\]

Consider lattice points $  (x_1,\ldots, x_n)$ 
and $  (y_1,\ldots, y_n)$ in the cube  $[0, h\ell]^n$.
Because $g > h\ell$,   the uniqueness of the $g$-adic 
representation of an integer implies 
\[
\varphi_g(x_1,\ldots, x_n) = \sum_{i=1}^n x_i g^{i-1} = \sum_{i=1}^n y_i g^{i-1} 
= \varphi_g(y_1,\ldots, y_n) 
\]
if and only if $x_i=y_i$ for all $i=1,\ldots, n$, and so the function $\varphi_g$ 
is one-to-one on the lattice point cube $ [0, h\ell]^n \cap \Z^n$.  
Because $\varphi_g$ is  linear and one-to-one, we have 
\beq          \label{lattice:sumset-1}
\varphi_g(hA)  = h\varphi_g(A) 
\eeq 
and 
\beq           \label{lattice:sumset-2}
\varphi_g \left( \widehat{hA} \right)  = \widehat{h\varphi_g(A) } 
\eeq
for every set $A$ in $[0, h\ell]^n \cap \Z^n$.  

If $t \in \mcr_{\Z^n} (h,k)$, then there is a finite set $A$ of lattice points in $\Z^n$ 
with $|A| = k$ and $|hA| = t$.    
Because sumset size is translation invariant, 
we can assume that the lattice points in $A$ have nonnegative integral coordinates, 
that is,  $A$ is a finite subset of  $\N_0^n$.  Choose integers $\ell$ and $g$ such that 
$A \subseteq [0,\ell]^n \subseteq \R^n$ and $g > h\ell$.   The sumset $hA$ 
is contained in  the cube $ [0, h\ell]^n$.  
Because  the function $\varphi_g$ is one-to-one and   linear on $ [0, h\ell]^n \cap \N_0^n$,  
we have 
\[ 
\varphi_g(A) \subseteq \N_0 \qqand 
|\varphi_g(A)| = |A| = k. 
\]
Identity~\eqref{lattice:sumset-1} implies  
\[
t = |hA| = |\varphi_g(hA)| = |h\varphi_g(A)| \in \mcr_{\Z} (h,k).
\]
It follows that $\mcr_{Z^n} (h,k) \subseteq \mcr_{\Z} (h,k)$ 
and so $\mcr_{Z} (h,k)= \mcr_{\Z^n} (h,k)$. 

Using identity~\eqref{lattice:sumset-2}, we apply   the same argument to restricted sumsets 
and obtain $\widehat{\mcr}_{\Z} (h,k) = \widehat{\mcr}_{\Z^n} (h,k)$. 
This completes the proof. 
\end{proof} 

The function $\varphi_g$ constructed in the proof of Theorem~\ref{lattice:theorem:Z} 
is a classical example of a Freiman isomorphism (Nathanson~\cite{nath96bb}).

Let $h \geq 2$ and $k \geq 2$.  
 Because the set  $\mcr_{\Z}(h,k)$  is finite, 
there is an integer $N$ such that,  for all $t \in \mcr_{\Z}(h,k)$, 
 there exists $A \subseteq [0,N]$ with $|A| = k$ and $|hA| = t$.  
Similarly, the set  $\widehat{\mcr}_{\Z} (h,k) $ is finite and 
there is an integer  $\widehat{N}$ such that,  for all $t \in \widehat{\mcr}_{\Z} (h,k) $, 
 there exists $A \subseteq [0,\widehat{N}]$ 
with $|A| = k$ and  $\left| \widehat{hA} \right| = t$.

There are the following questions (Nathanson~\cite{nath26q} and Gowers~\cite{gowe26}).

\bprob
Compute  integers $N$ and $\widehat{N}$  with the property that 
\beq          \label{lattice:N} 
\mcr_{\Z} (h,k)  = \left\{ \left| hA \right|: A \subseteq [0,N] \text{ and }  |A| = k \right\} 
\eeq 
and 
\beq          \label{lattice:N-hat} 
\widehat{\mcr}_{\Z}(h,k) = \left\{ \left| \widehat{hA} \right|: 
A \subseteq [0, \widehat{N}] \text{ and }  |A| = k \right\}. 
\eeq 
\eprob 

\bprob
Estimate or compute the smallest integers $N = N(h,k)$ 
and  $\widehat{N} = \widehat{N}(h,k)$ that satisfy~\eqref{lattice:N} 
and~\eqref{lattice:N-hat}. 
\eprob

By Theorem~\ref{lattice:theorem:Z}, we can compute $\mcr_{\Z} (h,k)$ 
and $\widehat{\mcr}_{\Z}(h,k) $ 
by constructing sets of lattice points and  computing $\mcr_{\Z^n} (h,k)$ 
and   $\widehat{\mcr}_{\Z^n}(h,k)$.  
It is natural to  look for upper and lower bounds 
for the diameters of finite sets of $k$ lattice points that suffice to compute 
the range of $h$-fold  sumset sizes and the range of  $h$-fold  restricted sumset sizes. 
There are different ways to define ``diameter.'' 
For simplicity, we choose the $\ell^{\infty}$-norm:  If $\mbx = (x_1,\ldots, x_n) \in \R^n$, then 
\[
\| \mbx \|_{\infty} = \|(x_1,\ldots, x_n)\|_{\infty} = \max\{|x_i|:i=1,\ldots, n\}. 
\]
The diameter of a set $A$ in $\R^n$ is 
\[
\diam(A) = \sup\{\| \mbx - \mby\|_{\infty} : \mbx, \mby \in A\}.
\]
If a set $A$ has diameter $d$, then, by translation, we can assume that $A \subseteq [0,d]^n$.  
The geometrical question is whether the freedom of working in $n$-dimensions 
enables the computation of the sumset size sets $\mcr_{\Z}(h,k)$ and $\widehat{\mcr}_{\Z^n}(h,k)$ 
by considering only sets with sufficiently small diameter that it sufices to examine 
fewer than $N(h,k)$  and $\widehat{N}(h,k)$ sets, respectively. 
One precise formulation is the following:

\bprob
Let $N(h,k)$ be the smallest positive integer such that 
\[
\mcr_{\Z} (h,k)  = \left\{ \left| hA \right|: A \subseteq [0,N(h,k)] \text{ and }  |A| = k \right\}. 
\]
and let $N_n(h,k)$ be the smallest positive integer such that 
\[
\mcr_{\Z} (h,k)  = \left\{ \left| hA \right|: A \subseteq [0,N_n(h,k)]^n \text{ and }  |A| = k \right\}. 
\]
Prove or disprove the following inequality:
\[
N_n(h,k)\leq N(h,k)^{1/n}.
\]
\eprob


\begin{thebibliography}{99} 

\bibitem{gowe26}
W. T. Gowers, A recent experience with ChatGPT 5.5 Pro, 
Gowers's Weblog, May 8, 2026

\bibitem{nath96bb} 
M.  B. Nathanson,
\emph{Additive Number Theory: Inverse Theorems and the Geometry of Sumsets},
Springer-Verlag, 1996. 

\bibitem{nath25a} 
M.  B. Nathanson, 
Problems in additive number theory, VI: Sizes of sumsets of finite sets,  
Acta Math.  Hungarica 176 (2025), 498--521.

\bibitem{nath26q} 
M.  B. Nathanson, 
Diversity, equity, and inclusion for problems  in additive number theory, 
arXiv:  2603.15556 
\end{thebibliography}
\end{document}